\documentclass[11pt]{article}

\usepackage{amsmath,amssymb,amsthm}
\usepackage{enumerate}
\usepackage[all]{xy}
\usepackage{xypic}

\def\tytul#1{\noindent{\bf\LARGE{#1}} \bigskip \thispagestyle{plain}}
\def\autor#1{\noindent{\large{#1}}\smallskip}
\def\wydzial#1{\smallskip \noindent \small #1}
\def\email#1{\smallskip \noindent \small E-mail address: {\sf #1}}

\theoremstyle{plain} \newtheorem{theorem}{Theorem}[section]
\theoremstyle{plain} \newtheorem{lemma}[theorem]{Lemma}
\theoremstyle{plain} \newtheorem{proposition}[theorem]{Proposition}
\theoremstyle{plain} \newtheorem{corollary}[theorem]{Corollary}
\theoremstyle{plain} \newtheorem{remark}[theorem]{Remark}
\theoremstyle{definition} \newtheorem{example}[theorem]{Example}
\theoremstyle{definition} \newtheorem{definition}[theorem]{Definition}
\theoremstyle{definition} \newtheorem{question}{Question}

\newcommand{\sierp}{\mathbb{S}}
\newcommand{\id}{\operatorname{id}}
\newcommand{\cC}{c$\mathcal{C}$}

\begin{document}

\tytul{On homotopy types of Alexandroff spaces}

\autor{Micha{\l} Kukie{\l}a}

\wydzial{Faculty of Mathematics and Computer Science,\\ Nicolaus Copernicus University,\\ ul. Chopina 12{\slash}18,\\ 87-100 Toru\'{n}, Poland}

\email{mckuk@mat.uni.torun.pl}

\begin{abstract}
We generalise some results of R. E. Stong concerning finite spaces to wider subclasses of Alexandroff spaces. These include theorems on function spaces, cores and homotopy type. In particular, we characterize pairs of spaces $X,Y$ such that the compact-open topology on $C(X,Y)$ is Alexandroff, introduce the classes of finite-paths and bounded-paths spaces and show that every bounded-paths space and every countable finite-paths space has a core as its strong deformation retract. Moreover, two bounded-paths or countable finite-paths spaces are homotopy equivalent if and only if their cores are homeomorphic. Some results are proved concerning cores and homotopy type of locally finite spaces and spaces of height 1. We also discuss a mistake found in an article of F.G. Arenas on Alexandroff spaces.

It is noted that some theorems of G. Minian and J. Barmak concerning the weak homotopy type of finite spaces and the results of R. E. Stong on finite H-spaces and maps from compact polyhedrons to finite spaces do hold for wider classes of Alexandroff spaces.

Since the category of $T_0$ Alexandroff spaces is equivalent to the category of posets, our results may lead to a deeper understanding of the notion of a core of an infinite poset.
\end{abstract}


%
%

\begin{section}{Introduction}
Alexandroff spaces, first introduced by P. Alexandroff \cite{alexandroff}, are topological spaces which have the property that intersection of any family of their open subsets is open. Finitie spaces are an important special case. A basic property of the category of Alexandroff spaces (and continuous maps) is that it is isomorphic to the category of preorders (and order-preserving maps), while the category of $T_0$ Alexandroff spaces is isomorphic to the category of posets.

From the point of view of general topology, Alexandroff spaces seem not to have drawn much interest since Alexandroff's work (there are, however, some articles on this topic, for example \cite{arenas} or \cite{mahdi}). On the other hand, they have gained much attention because of their use in digital topology, cf. \cite{kopperman}, \cite{melin}.

The situation is somewhat different if we turn to algebraic topology. Homotopy type and weak homotopy type of Alexandroff spaces were studied by several authors. For finite spaces, these were investigated for example by Stong in \cite{stong}, Minian and Barmak in \cite{barmak2}, \cite{barmak}. The latter two authors also developed the basics of simple homotopy theory for finite spaces in \cite{barmak1}. In \cite{mccord} McCord shows quite a surprising result -- for every polyhedron exists an Alexandroff space that is weak homotopy equivalent to it (and \textit{vice versa}). Also \cite{hardie} contains some results on this theme.

In this article we concentrate mostly on generalising the results of Stong \cite{stong}, which we shall now shortly describe. Stong's article contains results concerning homotopies of maps between finite spaces and homotopy type classification of finite spaces via their \textit{cores}. Cores are also known in order theory (under the same name) and were re-discovered for this discipline by Duffus and Rival \cite{duffus-rival}. From our point of view, cores are a special kind of deformation retracts. For a finite space, its core may be achieved by removing from the space, one by one, points known as beat points. (Our terminology comes from May's notes \cite{may}, it was also used by Minian and Barmak in \cite{barmak}. Stong called the points linear and colinear, while in order theory they are called irreducible points. We decided to use the name ``beat points'' mainly because it is the shortest one.) It turns out that two finite spaces are homotopy equivalent if and only if their cores are homeomorphic. An attempt to generalise these results to locally finite spaces was made by Arenas in \cite{arenas}. Unfortunately, his paper contains a mistake that makes the theory look a lot simpler than it really is. (Arenas states that $C(X,X)$ with the compact-open topology is Alexandroff for $X$ Alexandroff. We show that this is never true if $X$ is infinite.) In this article we take roughly the same route, but avoid falling into the same trap.

In Sections \ref{sec:preliminaries} and \ref{sec:spaces} we prepare the ground for further results. Basic facts concerning Alexandroff spaces are recalled and the classes of locally finite spaces, fp-spaces and bp-spaces are introduced. Section \ref{sec:homotopy} is a study of the compact-open topology on the space of continuous functions $C(X,Y)$ between two Alexandroff spaces. We estabilish a bijection between paths in $C(X,Y)$ and homotopies. Then we describe how to construct some of these paths and apply the results to show contractibility of the Khalimsky line (also known under the name ``infinite fence''). It is also here that the mistake in \cite{arenas} is discussed. In Section \ref{sec:cores} we investigate cores of spaces belonging to the classes introduced ind Section \ref{sec:spaces}. In particular, we show that every fp-space has a core and that two bp-spaces or countable fp-spaces are homotopy equivalent if and only if their cores are isomorphic. This may be seen as an another approach to the results of Section 3 of Farley's article \cite{farley}. The last section contains some theorems on $\gamma$-points (a concept introduced by Barmak and Minian in \cite{barmak} and related to the weak homotopy type of Alexandroff spaces), Alexandroff H-spaces and maps from compact polyhedrons into Alexandroff spaces. Some of the results presented in this article, especially in Section \ref{sec:cores}, may be known to order-theoretists, but they seem not to have been applied to studying homotopy types of infinite Alexandroff spaces.

The author believes the topological approach to cores of infinite posets may yield new insights into their, still not well understood, nature. On the other hand, expressing topological notions, such as homotopy, by the means of the language of order theory greatly simplifies their investigation. 
\end{section}

%
%

\begin{section}{Preliminaries}\label{sec:preliminaries}
We will now set up some terminology and recall a few basic properties of Alexandroff spaces. Some references for this section are \cite{mahdi}, \cite{may}, \cite{stong}, though the latter two articles are concerned with finite spaces only. Nevertheless, the following results may be proved for all Alexandroff spaces the same way it is done in the finite case. The reader is assumed to have some background in order theory, general and algebraic topology. 

A topological space $X$ is called an \textit{Alexandroff space} iff arbitrary intersections of sets open in $X$ are open. For every $x\in X$ the intersection of all open neighbourhoods of $x$ is then the minimal open neighbourhood of $x$. We shall denote it by $U_x$. The family $\{U_x\}_{x\in X}$ is a basis for the topology on $X$. (Existence of minimal open neighbourhoods of every point in a space is in fact equivalent to that space being Alexandroff.) By $\mathbf{Al}$ we denote the category of Alexandroff spaces and continuous maps and by $\mathbf{T_0 Al}$ its full subcategory of $T_0$ Alexandroff spaces.

A \textit{preorder} is a set equipped with a transitive, reflexive, binary relation. A \textit{partially ordered set} (a \textit{poset}) is a preorder such that the binary relation is also antisymmetric. Given two elements $p,q$ of a preorder $(P,\leq)$ by $p\sim q$ we will denote the fact that they are \textit{comparable}, that is $p\leq q$ or $q\leq p$. A poset $P$ is called \textit{chain-complete}, if every chain in $P$ has both a supremum and an infimum in $P$. For a poset $P$ by $\max (P)$ we denote the set of elements maximal in $P$ and by $\min (P)$ the set of elements minimal in $P$. For $x\in P$ we define $x\!\downarrow=\{y\in P:y\leq x\}$. Subset $A$ of a poset $P$ is called a \textit{down-set} if $a\!\downarrow\subseteq A$ for every $a\in A$. By $\mathbf{Preorder}$ we denote the category of preorders and order-preserving maps and by $\mathbf{Poset}$ its full subcategory of partially ordered sers. 

With every topological space $X$ we associate its \textit{specialization preorder}, that is a preorder $\mathcal{P}(X)=(X,\leq)$ such that $x\leq y$ iff $y\in\overline{\{x\}}$. This association is functorial. Continuous functions are mapped to increasing functions (with the same graph). If $X$ is $T_0$, then its specialization preorder is a partial order. If $X$ is $T_1$, then its specialization preorder is an antichain. 

With every preorder $(P,\leq)$ we associate an Alexandroff space $\mathcal{X}(P)$, with topology generated by the open basis $\{U_x=x\!\downarrow\}_{x\in X}$. Therefore, a set is open in $\mathcal{X}(P)$ iff it is a down-set in $(P,\leq)$. $\mathcal{X}$ is also a functor, increasing functions between preorders are continuous functions between the associated spaces. If $P$ is a partial order, then $\mathcal{X}(P)$ is $T_0$. 

If we restrict $\mathcal{P}$ to the category $\mathbf{Al}$, then the functors $\mathcal{X}:\mathbf{Preorder}\to\mathbf{Al}$ and $\mathcal{P}:\mathbf{Al}\to\mathbf{Preorder}$ are mutually inverse. The same property holds if we exchange $\mathbf{Preorder}$ with $\mathbf{Poset}$ and $\mathbf{Al}$ with $\mathbf{T_0 Al}$. Therefore, in this article we will not make a difference between an Alexandroff space and its associated preorder. The language used will come both from topology and order-theory. 

Starting from this point, all spaces are assumed to be $T_0$, unless stated otherwise. Since it may be shown the Kolmogorov quotient of an Alexandroff space is homotopy equivalent to that space (this is a consequence of the fact that Lemma \ref{finite_homotopies_lemma} works also for non-$T_0$ spaces), from our point of view little generality is lost, while the exposition gets clearer and shorter.

By $I$ we denote the unit interval $[0,1]\subseteq \mathbb{R}$ with the Euclidean topology. The following lemma is a well known result.

\begin{lemma}\label{path_lemma}
Let $X$ be an Alexandroff space and $x,y\in X$ be such that $x\leq y$. Then there exists a path $h:I\to X$ with $h(0)=x$ and $h(1)=y$. 
\end{lemma}
\begin{proof}
Let $h(t)=\begin{cases}x & \text{for } t<1\\ y & \text{for } t=1\end{cases}$. We easily check that $h$ is continuous.
\end{proof}

If $x_1,x_2,\ldots x_n\in X$ are such that $x_i\sim x_{i+1}$ for every $i=1,2\ldots,n-1$, then from the lemma it follows that there exists a path in $X$ from $x_1$ to $x_n$. In particular $x\!\downarrow=U_x$ is path connected and thus Alexandroff spaces are locally path connected. It is straightforward to check that for any $x\in X$ the set of $y\in X$ such that exists a finite sequence $x=x_1\sim x_2\sim\ldots\sim x_n=y$ is clopen. Therefore, it is the connected component (and the path component) of $x$.

A \textit{Khalimsky line} (or a \textit{two-way infinite fence}) is a space homeomorphic to $(\mathbb{Z},\sqsubseteq)$, where $n\sqsubseteq m$ iff either $n=m$ or $|n-m|=1$ and $m$ is even. A \textit{Khalimsky half-line} (or a \textit{one-way infinite fence}) is a space homeomorphic to one of the following two subpaces of $(\mathbb{Z},\sqsubseteq)$: $X_1=\{n:n\geq 0\}$, $X_2=\{n:n>0\}$, where $\geq, >$ refer to the standard order on $\mathbb{Z}$. Connected, finite subsets of the Khalimsky line are called \textit{finite fences}. Such spaces are important both in order theory and in digital topology (where they serve as an analogue of intervals of the real line, cf. \cite{kopperman}). 
\end{section}

%
%

\begin{section}{Some classes of Alexandroff spaces}\label{sec:spaces}
In this section we introduce a few subclasses of Alexandroff spaces and prove some of their properties. Each of these classes is, in a sense, close to the class of finite spaces. Keep in mind that all spaces are assumed $T_0$.

\begin{definition}
We call a sequence (finite or infinite) $(x_n)$ of elements of an Alexandroff space $X$ an \textit{s-path} if $x_i\not=x_j$ for $i\not=j$ and $x_{i-1}\sim x_{i}$ for all $i>0$. Given a finite s-path $s=(x_0,x_1,\ldots,x_m)$ we say $m$ is the \textit{length} of $s$ and call $s$ an \textit{s-path from $x_0$ to $x_m$}.
\end{definition}

\begin{definition} We say an Alexandroff space $X$ is:
\begin{enumerate}
\item a \textit{finite-chains space}, if every chain in $X$ is finite,
\item a \textit{locally finite space}, if for every $x\in X$ the set $\{y\in X: y\sim x\}$ is finite,
\item a \textit{finite-paths space} (\textit{fp-space}), if every s-path of elements of $X$ is finite,
\item a \textit{bounded-paths space} (\textit{bp-space}), if exists $n\in\mathbb{N}$ such that every s-path of elements of $X$ has less than $n$ elements.
\end{enumerate}
\end{definition}

Note that different notions of local finiteness are in use, which do not coincide and should not be confused with our definition. For example, sometimes it is only required that $\{y:y\leq x\}$ is finite for every $x\in X$. Also, the name is sometimes used for posets with the property that all intervals $[x,y]=\{z:x\leq z\leq y\}$ are finite. 

Obviously, bp-spaces form a (strict) subclass of fp-spaces and both fp-spaces and locally finite spaces are (strict) subclasses of finite-chains spaces. Moreover, from Proposition \ref{loc_fin_spanning} below it follows that connected spaces which are both locally finite and fp-spaces must be finite.

\begin{definition}
Recall from \cite{duffus} we call a subset $A$ of a poset $X$ \textit{spanning} if $A\subseteq \min(X)\cup\max(X)$. By $d_X(x,y)$ we mean the minimal length of an s-path in $X$ from $x$ to $y$ (or $\infty$ if such an s-path does not exist). $A$ is called \textit{isometric} if $d_A(x,y)=d_X(x,y)$ for all $x,y\in A$. By $B(x,n)$ we denote the set $\{y\in X:d_X(x,y)\leq n\}$. 
\end{definition}

\begin{proposition}\label{loc_fin_spanning}
A connected, locally finite space is infinite if and only if it contains an isometric, spanning subset isomorphic to the Khalimsky half-line.
\end{proposition}
\begin{proof}
The "if" part is trivial.

For the "only if" part, suppose $X$ is a connected, locally finite space and fix a point $x_0\in \max(X)\cup\min(X)$. Since $X$ is connected, for every $y\in X$ exists an isometric fence $x_0=a_0^y\sim a_1^y \sim \ldots \sim a_{n(y)}^y=y$ with $n(y)=d_X(x_0,y)$. Without lack of generality we may assume $a_m^y\in\max(X)\cup\min(X)$ for all $m<n(y)$. Let $f:X\smallsetminus\{x_0\}\to X$ be given by $f(y)=a^y_{n(y)-1}$. For every $y\in X$ the set $f^{-1}(y)$ is finite, as a subset of the finite set $B(y,1)$. Moreover, it is easy to see that the fence \[F(y)=\{x_0=f^{n(y)}(y),f^{n(y)-1}(y),\ldots,y\}\label{seq1}\] is isometric. Since $X$ is locally finite, $B(x_0,n)$ is finite for every $n>0$. Thus $X\not\subseteq B(x_0,n)$ for any $n$. Therefore, for every $n\in\mathbb{N}$ exists an $y\in X$ with $d_X(x_0,y)=n$. By the K{\"o}nig's lemma (cf. Chapter III, §5 of \cite{kuratowski}) there exists an infinite sequence $S=(a_0=x_0,a_1,a_2,\ldots)$ with $a_{n-1}=f(a_n)$ for all $n\geq 1$. It is isometric, which follows from the fact that all fences $F(y)$ are isometric. Moreover, we have assumed $f(y)\in\max(X)\cup\min(X)$ for all $y\in X$, thus the one-way infinite fence $S$ is spanning.
\end{proof}

We get a simple corollary.

\begin{corollary}\label{loc_fin_non_fpp}
An infinite, locally finite space does not have the fixed point property.
\end{corollary}
\begin{proof}
Let $X$ be an infinite, locally finite space. If $X$ is not connected, then $X$ does not have the fixed point property. If $X$ is connected, the result follows from the above Proposition \ref{loc_fin_spanning} and Theorem 3 of \cite{duffus}.
\end{proof}

\begin{definition}
We say an Alexandroff space $X$ satisfies:
\begin{enumerate}
\item the \textit{descending chain condition} (DCC), if it contains no strictly decreasing infinite sequence,
\item the \textit{ascending chain condition} (ACC), if it contains no strictly increasing infinite sequence.
\end{enumerate}
\end{definition}

Proofs of the following two lemmas are standard.
\begin{lemma}\label{fin_chain}
An Alexandroff space is a finite-chains space iff it satisfies both ACC and DCC. 
\end{lemma}

\begin{lemma}\label{max_is_acc}
An Alexandroff space $X$ satisfies ACC iff for every subset $A$ of $X$ the set $\max(A)$ is not empty. (By duality, the same result holds for spaces satisfying DCC and sets $\min(A)$.)
\end{lemma}

\end{section}

%
%

\begin{section}{Function spaces}\label{sec:homotopy}
If $X, Y$ are topological spaces, by $C(X,Y)$ we denote the set of continuous functions from $X$ to $Y$ equipped with the compact-open topology.

\begin{lemma}[cf. \cite{fox}]
If $X, T$ are topological spaces which satisfy the first countability axiom and $Y$ is an arbitrary topological space, then continuity of $h:X\times T\to Y$ is equivalent to continuity of $h^*:T\to C(X,Y)$, where $h^*(t)(x)=h(x,t)$.
\end{lemma}

The above lemma holds for $X$ being an Alexandroff space (each point has a neighbourhood basis consisting of one element -- the minimal neighbourhood, so the space is first-countable) and $T=I$, the unit interval. Therefore, we have the following corollary.

\begin{corollary}\label{path_is_homotopy}
Let $X$ be an Alexandroff space, $Y$ an arbitrary topological space. Maps $f,g:X\to Y$ are homotopic if and only if they belong to the same path component of $C(X,Y)$. 
\end{corollary}

On the space $C(X,Y)$, where $Y$ is Alexandroff and $X$ an arbitrary topological space, we define an order $f\leq g$ iff $f(x)\leq g(x)$ for all $x\in X$. The following proposition may be proved as Proposition 9 is in \cite{stong}. 
\begin{proposition}
Let $X$ be a topological space and $Y$ an Alexandroff space. The intersection of all open sets in $C(X,Y)$ containing the map $f$ is equal to $f\!\downarrow=\{g\in C(X,Y):g\leq f\}$.
\end{proposition}
Therefore, $g\leq f$ implies $f\in \overline{\{g\}}$, so the order on $C(X,Y)$ defined above is the specialization order of $C(X,Y)$. In Theorem 3.1 of \cite{arenas} it is stated that $C(X,Y)$ is an Alexandroff space, which in general is not true (intersection of all open sets containing a map doesn't have to be open). We will show this in Theorem \ref{arenas_mistake}. As a consequence of the mistake, several other untrue statements are made in \cite{arenas}; we will point them out further in the text. 

By $[A,B]$, where $A\subseteq X$, $B\subseteq Y$, we will denote the set $\{f\in C(X,Y):f(A)\subseteq B\}$. For brevity, by $[x,U]$ we denote $[\{x\},U]$.

\begin{proposition}
If $X$ is an Alexandroff space, then $X$ is compact if and only if $\max(X)$ is finite and for every $x\in X$ exists an $y\in\max (X)$ such that $y\geq x$.
\end{proposition}
\begin{proof}
For the "if" part, let $\mathcal{U}$ be an open cover of $X$. For every $y\in\max(X)$ choose an open set $U(y)\in\mathcal{U}$ such that $y\in U(y)$. If for every $x\in X$ exists an $y\in\max (X)$ such that $y\geq x$, then (since open sets are down-sets) $\{U(y):y\in\max(X)\}$ is a subcover of $\mathcal{U}$. If $\max(X)$ is finite, then the subcover is finite.

To prove the reverse implication, let $\mathcal{U}=\{x\!\downarrow:x\in X\}$. First suppose that $\max(X)$ is infinite. Then $\mathcal{U}$ does not have a finite subcover, since the only set in $\mathcal{U}$ containing an $x\in\max(X)$ is $x\!\downarrow$, so every subcover of $\mathcal{U}$ contains an infinite subset $\{x\!\downarrow:x\in\max(X)\}$. Now suppose there exists an $y\in X$ such that for every $x\geq y$ exists a $z>x$. If $\mathcal{V}=\{x_1\!\downarrow, x_2\!\downarrow,\ldots, x_n\!\downarrow\}$ is a finite subcover of $\mathcal{U}$, then we can choose an $x_i$ such that $x_i\!\downarrow\in\mathcal{V}$, $x_i\geq y$ and it is maximal among elements having this property. But then exists a $z>x_i$. By choice of $x_i$, $z\not\in\bigcup\mathcal{V}$. This is a contradiction with $\mathcal{V}$ being a cover.
\end{proof}

Now, let's use our characterization of compact sets to examine the compact-open topology on $C(X,Y)$ for $X,Y$ Alexandroff. 

\begin{lemma}\label{subbasis_lemma}
Let $X, Y$ be Alexandroff spaces, $K\subseteq X$ be compact, $U\subseteq Y$ be open. Then $[K,U]=[\max (K), U]$.
\end{lemma}
\begin{proof}
$\max (K)\subseteq K$, so $[K, U]\subseteq [\max (K),U]$. For the reverse inclusion, notice that $f(K)\subseteq U$ iff $f(x)\in U$ for all $x\in K$. But for every $x\in X$ exists an $x'\in\max (K)$ with $x'\geq x$. Because $f$ is order preserving, $f(x')\geq f(x)$, and since $U$ is a down set, $f(x)\in U$ if $f(x')\in U$. Therefore, $[\max (K), U]\subseteq [K,U]$.
\end{proof}

\begin{corollary}\label{product_corollary}
Let $X, Y$ be Alexandroff spaces. The family \mbox{$\{[x,y\!\downarrow]\}_{x\in X,y\in Y}$} is a subbasis for the compact-open topology on $C(X,Y)$. Therefore, the compact-open topology on $C(X,Y)$ coincides with the topology induced from the product space $\prod_{x\in X} Y$.
\end{corollary}
\begin{proof}
From Lemma \ref{subbasis_lemma} we know that sets of the form $[\{x_1,\ldots,x_n\},U]$, where $x_1,\ldots,x_n\in X$ and $U\subseteq Y$ is open, form a subbasis for $C(X,Y)$. Because $[\{x_1,\ldots,x_n\},U]=\bigcap_{i=1}^{n}[x_i,U]$, sets of the form $[x,U]$ also constitute a subbasis. Now, for a basic set $\bigcap_{i=1}^{n}[x_i,U_i]$ we have:
\[\bigcap_{i=1}^{n}[x_i,U_i]=\bigcap_{i=1}^{n}\bigcup_{u\in U_i}[x_i,u\!\downarrow]=\bigcup_{u_1\in U_1}\ldots\bigcup_{u_n\in U_n}\bigcap_{i=1}^{n}[x_i,u_i\!\downarrow],\]
so $\{[x,y\!\downarrow]\}_{x\in X,y\in Y}$ is a subbasis for the compact-open topology on $C(X,Y)$.

The second statement of the corollary follows from the above result and the definition of the product topology.
\end{proof}

\begin{corollary}
Let $X$ be a finite space, $Y$ an Alexandroff space. Then the compact-open topology on $C(X,Y)$ is Alexandroff. 
\end{corollary}
\begin{proof}
This follows from Corollary \ref{product_corollary} and the (easily-checked) fact that finite products and subsets of Alexandroff spaces are Alexandroff.
\end{proof}

\begin{corollary}
Let $X,Y$ be Alexandroff spaces. Then the compact-open topology on $C(X,Y)$ is weaker than the Alexandroff topology induced by order on $C(X,Y)$.
\end{corollary}
\begin{proof}
Since $[K,U]$ are down-sets and both unions and intersections of down-sets are again down-sets, open sets in $C(X,Y)$ are all down-sets and therefore are open in the Alexandroff topology.
(This corollary may also be deduced from the general fact that Alexandroff topology is the strongest topology inducing given preorder.)
\end{proof}

\begin{lemma}\label{hereditarily_compact}
An Alexandroff space $X$ is hereditarily compact if and only if it contains no infinite antichains and satisfies the ACC.
\end{lemma}
\begin{proof}
If $X$ is hereditarily compact, then it contains no infinite antichains and no infinite, strictly increasing sequences, since these sets are not compact.

To prove the reverse, notice that if $X$ contains no infinite antichains and satisfies ACC, then every subset of $X$ also has this property. Therefore, it suffices to show that a space $X$ satisfying ACC and containing no infinite antichains is compact. By Lemma \ref{max_is_acc}, $\max (x\!\uparrow)\not=\emptyset$ for every $x\in X$. But $\max (x\!\uparrow)\subseteq \max (X)$, so every element of $X$ is under some maximal element. If the space also contains no infinite antichains, then the set of its maximal points is finite (because it is an antichain) and the space is compact.
\end{proof}

\begin{theorem}\label{arenas_mistake}
Let $X$ be an infinite Alexandroff space, $Y$ be an arbitrary Alexandroff space with $|Y|\geq 2$. Then the following statements are true.
\begin{enumerate}
\item\label{pt1} If $Y$ is discrete, then $C(X,Y)$ is Alexandroff if and only if $X$ has finitely many connected components.
\item\label{pt2} If $Y$ is not discrete, then $C(X,Y)$ is Alexandroff if and only if $X$ is hereditarily compact and for every $f\in C(X,Y)$ the image $f(X)$ is finite. 
\end{enumerate}
\end{theorem}
\begin{proof}
Since we know that the compact-open topology on $C(X,Y)$ is always weaker than the Alexandroff topology, we will only investigate for what spaces it is also stronger, and that is the case if and only if $f\!\downarrow$ is open for every $f\in C(X,Y)$.

For the proof of \ref{pt1}, suppose $Y$ is discrete and $f\in C(X,Y)$. Since $f$ is continuous, it maps every connected component of $X$ to a point. If $X$ has finitely many connected components, say $S_1,\ldots,S_n$, for every component $S_i$ we may choose an $s_i\in S_i$. Then $f\!\downarrow=\bigcap_{i=1}^{n}[s_i,f(s_i)\!\downarrow]$, which is open in $C(X,Y)$. Now suppose $X$ has infinitely many connected components. We will show that any basic set $\bigcap_{i=1}^{n}[x_i,y_i\!\downarrow]$ contains a map that is not in $f\!\downarrow$. For every set of this form a connected component $S$ of $X$ exists such that $\{x_1,\ldots,x_n\}\cap S=\emptyset$. Since $|Y|\geq 2$, there is an $y\in Y$ such that $S$ is not mapped to the point $y$. Let $f_0(x)=\begin{cases}f(x) & \text{for } x\not\in S\\ y & \text{for } x\in S\end{cases}$. Then $f\!\downarrow\not\ni f_0\in \bigcap_{i=1}^{n}[x_i,y_i\!\downarrow]$.

To prove \ref{pt2}, first assume $X$ is not hereditarily compact. Notice that if $Y$ is not discrete, then it contains a copy of the Sierpi{\'n}ski space (that is the space $\sierp=\{0,1\}$ with topology $\{\emptyset,\sierp,\{0\}\}$) as its subspace.  Now, since $X$ is not hereditarily compact, a subspace $A\subseteq X$ exists that is not compact. Let $A\!\downarrow=\bigcup_{a\in A}(a\!\downarrow)$. Define $f:X\to \sierp\hookrightarrow Y$ by $f(A\!\downarrow)=0$ and $f(X\smallsetminus (A\!\downarrow))=1$. Then $f^{-1}(\{0\})=A\!\downarrow$, which means the preimage of the only non-trivial open set in $\sierp$ is a down-set in $X$. Therefore $f$ is continuous. For a basic neighbourhood $\bigcap_{i=1}^{n}[x_i,y_i\!\downarrow]$ of $f$, $y_i\geq 1$ for all $x_i\not\in A\!\downarrow$. Since $A$ is not compact, an $a\in A$ exists such that $a\not\leq x_j$ for all $x_j\in A\!\downarrow$. Now define $f_0\in \bigcap_{i=1}^{n}[x_i,y_i\!\downarrow]$ by $f_0(x)=\begin{cases}f(x) & \text{for } x\not\geq a\\ 1 & \text{for } x\geq a\end{cases}$. Since $(A\!\downarrow)\smallsetminus (a\!\uparrow)$ is a down set, $f_0$ is continuous. But $f_0\not\leq f$.

Suppose now that $X$ is hereditarily compact. If $f(X)$ is finite for an $f\in C(X,Y)$, then $f\!\downarrow=\bigcap_{y\in f(X)}[f^{-1}(y),y\!\downarrow]$ is an open set. Therefore, if $f(X)$ is finite for every $f\in C(X,Y)$, we are finished. Suppose a $g\in C(X,Y)$ exists such that $g(X)$ is infinite. $g(X)$ is hereditarily compact as a continuous image of a hereditarily compact space. Therefore $g(X)$ satisfies ACC and has no infinite antichains. It is a well-known fact that every infinite poset contains either an infinite chain or an infinite antichain (see Corollary 2.5.10 of \cite{schroder}). $g(X)$ is infinite, so it cannot satisfy DCC (otherwise, it would be a finite-chains space without infinite antichains, thus finite), which means $g(X)$ contains a descending chain $c_0>c_1>c_2>\ldots$. For every $n\in\mathbb{N}$ choose an $a_n\in g^{-1}(c_n)$. Then $\{a_n:n\in\mathbb{N}\}$ is an infinite, hereditarily compact subset of $X$ and, by the same reasoning as for $g(X)$, it must contain an infinite descending chain $a_{k_0}>a_{k_1}>a_{k_2}>\ldots$. We will now define a map $f\in C(X,Y)$. For $x\in X\smallsetminus (a_{k_0}\!\downarrow)$ let $f(x)=g(a_{k_0})$, for $x\in (a_{k_n}\!\downarrow)\smallsetminus (a_{k_{n+1}}\!\downarrow)$, where $n=0,1,\ldots$, let $f(x)=g(a_{k_n})$, and for $x\in\bigcap_{n=1}^{\infty}a_{k_n}\!\downarrow$ let $f(x)=g(x)$. From the construction it is easy to see that $f$ is order-preserving. Now take a basic neighbourhood $\bigcap_{i=1}^{n}[x_i,y_i\!\downarrow]$ of $f$. 
If there exists an $x_i\in (a_{k_0}\!\downarrow) \smallsetminus \bigcap_{n=1}^{\infty}(a_{k_n}\!\downarrow)$, then let $N$ be the greatest number $l$ such that there exists an $x_j\in (a_{k_l}\!\downarrow) \smallsetminus \bigcap_{n=1}^{\infty}(a_{k_n}\!\downarrow)$. Otherwise, let $N=k_0$. Define $f_0:X\to Y$ by 
\[f_0(x)=\begin{cases}g(a_{k_N}) & \text{for } x\in (a_{k_{N+1}}\!\downarrow)\smallsetminus\bigcap_{n=1}^{\infty}(a_{k_n}\!\downarrow) \\ f(x) & \text{otherwise}\end{cases}.\]
Then $f_0\in \bigcap_{i=1}^{n}[x_i,y_i\!\downarrow]$ is order-preserving and $f_0>f$. 
\end{proof}

\begin{corollary}
For $X$ an infinite Alexandroff space the space $C(X,X)$ is not Alexandroff.
\end{corollary}
\begin{proof}
Suppose $X$ is an infinite Alexandroff space and $C(X,X)$ is Alexandroff. If $X$ was discrete, then by part 1 of Theorem \ref{arenas_mistake}, $X$ would have only finitely many components and therefore it would be finite. Thus $X$ is not discrete. But by part 2 of the theorem, it means that every map in $C(X,X)$ has a finite image. In particular $\operatorname{id}_X(X)=X$ is finite. A contradiction.
\end{proof}

\begin{corollary}
If $X$ is a hereditarily compact Alexandroff space and $Y$ is an Alexandroff space satisfying the DCC, then $C(X,Y)$ is Alexandroff.
\end{corollary}
\begin{proof}
For every $f\in C(X,Y)$ the image $f(X)$ is a hereditarily compact DCC space, thus it contains no infinite chains and no infinite antichains, and, by Corollary 2.5.10 of \cite{schroder}, is finite.
\end{proof}

In Corollary 3.3 of \cite{arenas} it is stated that homotopy classes of maps from $X$ to $Y$ ($X,Y$ being Alexandroff) coincide with connected components of $C(X,Y)$. That would follow if $C(X,Y)$ was locally path connected, in particular if $C(X,Y)$ was Alexandroff. But, for example, if $X$ is an infinite discrete space, then $C(X,X)$ is not locally path connected, so we cannot prove the statement this way. On the other hand, in the given counterexample path components and connected components of $C(X,X)$ are trivially the same. Therefore, we pose the following question.
\begin{question}
For Alexandroff spaces $X,Y$, do the path components of $C(X,Y)$ coincide with the connected components of $C(X,Y)$? If it is not so in the general, what conditions on $X,Y$ do guarantee such a coincidence?
\end{question}

Using some of the above results on $C(X,Y)$ we shall now describe a class of homotopies between functions in $C(X,Y)$ (or rather a class of paths in $C(X,Y)$, which is equivalent, since Corollary \ref{path_is_homotopy} holds). We begin with the following lemma.

\begin{lemma}\label{finite_homotopies_lemma}
Let $X,Y$ be Alexandroff spaces. If $f,g\in C(X,Y)$ are such that $f(x)\sim g(x)$ for all $x\in X$, then $f$ is homotopic to $g$ by a homotopy that is constant on the set $\{x\in X:f(x)=g(x)\}$. 
\end{lemma}
\begin{proof}
If $f(x)\sim g(x)$ for all $x\in X$, then $f\geq f'\leq g$, where the map $f'(x)=\min(f(x),g(x))$ is order-preserving. Therefore, it is sufficient to prove the lemma for $f\leq g$.

If $f\leq g$, then, by Lemma \ref{path_lemma}, there is a path $h:I\to C(X,Y)^{Al}$ given by $h(t)=\begin{cases}f & \text{for } t<1\\ g& \text{for } t=1\end{cases},$ where $C(X,Y)^{Al}$ denotes the space $C(X,Y)$ with the Alexandroff topology. Since the compact-open topology is weaker than the Alexandroff topology, the map $h:I\to C(X,Y)$ is also continuous. Therefore, by Corollary \ref{path_is_homotopy}, it induces a homotopy.
\end{proof}

If we can find $f_0,f_1,\ldots,f_n\in C(X,Y)$ such that $f_i\sim f_{i+1}$ for $i=0,\ldots,n-1$, then $f_0\simeq f_n$. If $C(X,Y)$ were Alexandroff, all homotopies could be expressed in such a way. In general, however, this is not the case (see Example \ref{khalimsky_contractible}). 

\begin{theorem}\label{infinite_homotopies}
Let $X, Y$ be Alexandroff spaces. If there exists an infinite sequence of functions $f_1, f_2, \ldots\in C(X,Y)$ such that $f_i(x)\sim f_{i+1}(x)$ for all $x\in X, i\in\mathbb{N}$, and a function $f_{\omega}\in C(X,Y)$ such that for every $x\in X$ exists $N\in\mathbb{N}$ with the property that $f_n(x)\leq f_{\omega}(x)$ for all $n\geq N$, then $f_0$ is homotopic to $f_{\omega}$.
\end{theorem}
\begin{proof}
We may assume $f_1\sim f_2\sim \ldots$. Otherwise, the sequence may be modified like in the proof of Lemma \ref{finite_homotopies_lemma} to fulfil this condition.

For $i=1,2,\ldots$ let $h_i:I\to C(X,Y)$ denote the path from $f_i$ to $f_{i+1}$ such as constructed in Lemma \ref{finite_homotopies_lemma}. Explicitly, $h$ is given by 
\[h_i(t)(x)=\begin{cases}f_i(x) & \text{for } t<1\\ f_{i+1}(x) & \text{for } t=1\end{cases}\] if $f_i\leq f_{i+1}$ and \[h_i(t)(x)=\begin{cases}f_{i}(x) & \text{for } t=0\\ f_{i+1}(x) & \text{for } t>0\end{cases}\] if $f_i\geq f_{i+1}$.

We shall define a path $H:I\to C(X,Y)$ from $f_1$ to $f$. For $n=1,2\ldots$, $t\in[1-\frac{1}{n},1-\frac{1}{n+1}]$ let $H(t)=h_n(n(n+1)(t-\frac{n-1}{n}))$ and let $H(1)=f_{\omega}$. Obviously, $H$ is well-defined. To complete the proof, it remains to show that $H$ is continuous. It's easy to see that $H$ is continuous at every $t\in [0,1)$ as a ``composition'' of paths. Therefore, we only have to show continuity at 1. Let $[x,U]$ be any subbasic neighbourhood of $H(1)=f_{\omega}$. That means $f_{\omega}(x)\in U$. We know an $N\in\mathbb{N}$ exists such that $f_n(x)\leq f_{\omega}(x)$ for all $n\geq N$, and thus $f_n(x)\in U$ for $n\geq N$ and $H((1-\frac{1}{N},1])\subseteq [x,U]$, which completes the proof. 
\end{proof}

As a simple application of Theorem \ref{infinite_homotopies} we will now characterize a class of contractible spaces (including the Khalimsky line and half-line).

\begin{theorem}\label{one_path_contractible}
Let $X$ be a finite-chains space and $x\in X$ a point such that for every $y\in X$ exists one and only one s-path \[S_y=(x=s_y^0,s_y^1,\ldots,s_y^{n(y)}=y)\] with the property that one of the elements $s_y^i, s_y^{i+1}$ is a cover of the other for all $0\leq i< n(y)$. Then $x$ is a strong deformation retract of $X$.
\end{theorem}
\begin{proof}
Let $f_k:X\to X$ for $0\leq k <\omega$ be given by \[f_k(y)=\begin{cases}s^y_k & \text{for } k\leq n(y)\\ y & \text{for } k> n(y)\end{cases}.\] Obviously, $f_k(y)\sim f_{k+1}(y)$ for all $y\in X$.

We need to show the maps $f_k$ are order preserving. Suppose $y\leq y'\in X$. There is a finite chain $C$ with $y=\min(C)$ and $y'=\max(C)$. If $x\in C$, then it is easy to see $f_k(y)\leq f_k(y')$ for all $0\leq k$. Suppose $x\not\in C$. One of the s-paths $S_y$, $S_{y'}$ forms the first $\min(n(y),n(y'))$ elements of the other, and the further elements of the longer path come from $C$. Otherwise, we could define another s-paths $\overline{S}_y, \overline{S}_{y'}$ that would have the above property. Suppose $n(y)\leq n(y')$. (In the other case we proceed analogously.) If $k\leq n(y)$, then $f_k(y)=f_k(y')$. If $k\geq n(y)$, then $f_k(y)=y$ and $f_k(y')\in C$, thus $f_k(y)\leq f_k(y')$. 

Application of Theorem \ref{infinite_homotopies} gives a homotopy from $f_0=x$ to $f_\omega=\id_X$. 
\end{proof}

\begin{corollary}\label{when_height_one_is_contractible}
Let $X$ be an Alexandroff space of height 1. Then the following statements are equivalent.
\begin{enumerate}
\item\label{st1} $X$ has a point as its strong deformation retract.
\item\label{st2} $X$ is path connected and the fundamental group of $X$ is trivial.
\item\label{st3} $X$ is path connected and the first homology group of $X$ is trivial.
\item\label{st4} $X$ is connected and it contains no crowns.
\end{enumerate}
\end{corollary}
\begin{proof}
\ref{st1} implies \ref{st2} and \ref{st2} implies \ref{st3} for every topological space.

Lemma 4.4.4 of \cite{schroder} says that a crown of minimal cardinality in a poset of height 1 is its retract. Since the first singular homology group of a crown is non-trivial (cf. \cite{barmak} or \cite{mccord}), \ref{st3} implies \ref{st4}.

An Alexandroff space of height 1 satisfies the assumptions of Theorem \ref{one_path_contractible} if and only if it is connected and contains no crowns. Therefore, \ref{st4} implies \ref{st1}. 
\end{proof}

\begin{example}\label{khalimsky_contractible}
Let $X=\{n:n\geq 0\}$ be the Khalimsky half-line, as described in Section \ref{sec:preliminaries}. By Corollary \ref{when_height_one_is_contractible}, $X$ is contractible. (The same applies to the Khalimsky line.)

On the other hand, it is easy to see that the only continuous function $f^0:X\to X$ with $f^0\sim \id_X$ is given by $f^0(0)=1$ and $f^0(n)=n$ if $n>0$. Now, the only function $f^1\sim f^0$ other than identity is given by $f^1(0)=f^1(1)=2$ and $f^1(n)=n$ if $n>1$. Continuing the reasoning we get an infinite sequence of comparable functions, none of which is constant. Therefore, there is no finite sequence $\id_X\sim f^0\sim f^1\sim\ldots f^n$ with $f^n$ a constant map.
\end{example}

We shall prove a generalisation of Theorem \ref{infinite_homotopies}, which will be useful in the next section. 

\begin{theorem}\label{countable_homotopies}
Let $X, Y$ be Alexandroff spaces and $\{f_{\alpha}:X\to Y\}_{\alpha\leq\gamma}$, where $\gamma$ is a countable ordinal, be a family of continuous maps such that:
\begin{enumerate}
\item if $\alpha=\beta+1$, then $f_\alpha\sim f_\beta$;
\item if $\alpha$ is a limit ordinal, then for every $x\in X$ exists $\beta^{\alpha}_{x}<\alpha$ such that $f_{\beta}(x)\leq f_{\alpha}(x)$ for all $\beta^{\alpha}_{x}\leq\beta\leq\alpha$.
\end{enumerate}
Then $f_0$ is homotopic to $f_{\gamma}$.
\end{theorem}
\begin{proof}
It is a well-known fact that any countable total order is isomorphic to a suborder of the rational numbers (cf. Chapter VI, §3, Theorem 3 of \cite{kuratowski}), and thus of the real unit interval. Therefore, we will identify $\{\alpha:\alpha\leq\gamma\}$ with a subset of $I=[0,1]$. Moreover, we may assume that $\gamma=1$ and $0=0$ (where the first $0$ in this equality indicates the ordinal number and the second $0$ is the real number).

We now define $H:I\to C(X,Y)$. For every ordinal $\alpha<\gamma$ let the restriction $H|_{[\alpha,\alpha+1]}:[\alpha,\alpha+1]\to C(X,Y)$ be constructed from the path from $f_{\alpha}$ to $f_{\alpha+1}$ such as defined in Lemma \ref{finite_homotopies_lemma}, in the same way it is done in the proof of Theorem \ref{infinite_homotopies}. Let $H(\gamma)=f_{\gamma}$.

It now suffices to show that $H$ is continuous. This may be done by a similar argument as in the proof of Theorem \ref{infinite_homotopies}, but applied to every limit ordinal $\alpha\leq\gamma$.
\end{proof}

\end{section}

%
%

\begin{section}{Cores}\label{sec:cores}
We shall now discuss cores of Alexandroff spaces, which from our viewpoint are special kinds of deformation retracts.

\begin{definition}
Let $X$ (or: $(X,p)$) be an Alexandroff space (with the distinguished point $p$). A point $x\in X$ ($x\in X\smallsetminus\{p\}$) is called an \textit{up-beat point} if the set $x\!\uparrow\smallsetminus\{x\}$ has a smallest element $u_x$. Dually, $x$ is a \textit{down-beat point} if $x\!\downarrow\smallsetminus\{x\}$ has a largest element $d_x$. A point that is either an up-beat point or a down-beat point is called a \textit{beat point}.
\end{definition}

\begin{definition}
Let $X$ be an Alexandroff space (with the distinguished point $p$). A retraction $r:X\to A\subseteq X$ ($r:(X,p)\to (A,p)\subseteq (X,p)$) is called:
\begin{enumerate}
\item a \textit{comparative retraction}, if $r(x)\sim x$ for every $x\in X$,
\item an \textit{up-retraction}, if $r(x)\geq x$ for every $x\in X$,
\item a \textit{down-retraction}, if $r(x) \leq x$ for every $x\in X$,
\item a \textit{retraction removing a beat point}, if exists an $x\in X$ being an up-beat point under some $u_x$ or a down-beat point over some $d_x$ and such that $r(x)=u_x$ or $r(x)=d_x$, and $r(y)=y$ for all $y\not=x$.
\end{enumerate}
\end{definition}

\begin{remark}\label{retractions_remark}
Retractions removing a beat point, up- and down-retractions are all comparative. Moreover, every comparative retraction may be written as a composition of an up-retraction and a down-retraction: if $r:X\to A$ is a comparative retraction, then $r=r_d\circ r_u$, where 
\[r_u(x)=\begin{cases}r(x) & \text{if } r(x)\geq x\\ x & \text{if } r(x)\leq x\end{cases}\]
 and 
\[r_d(x)=\begin{cases}r(x) & \text{if } r(x)\leq x\\ x & \text{if } r(x)\geq x\end{cases}.\]
It is straightforward to check $r_u$ and $r_d$ are order-preserving.

Moreover, by Lemma \ref{finite_homotopies_lemma}, $i\circ r$, where $i:A\hookrightarrow X$, is homotopic to $id_X$ by a homotopy keeping fixed the set $A$, and this means $r$ is a strong deformation retraction. 
\end{remark}

\begin{definition}
In the category of Alexandroff spaces (with distinguished points) let $\mathcal{C}$ denote the class of all comparative retractions, $\mathcal{U}$ and $\mathcal{D}$ classes of, respectively, up- and down-retractions and $\mathcal{I}$ the class of retractions removing a beat point.
\end{definition}

\begin{definition}
Let $\mathcal{R}$ be a class of retractions in the category of Alexandroff spaces (with distinguished points). A nonempty Alexandroff space $X$ ($(X,p)$) is called an \textit{$\mathcal{R}$-core} if there is no retraction $r:X\to r(X)$ ($r:(X,p)\to (r(X),p)$) in $\mathcal{R}$ other than identity.
\end{definition}

By Remark \ref{retractions_remark}, a space being a $\mathcal{C}$-core is equivalent to that space being an $\mathcal{U}\cup\mathcal{D}$-core. Also, every $\mathcal{C}$-core is an $\mathcal{I}$-core. In general, the reverse implication does not hold, as we can see in the following example. 

\begin{example}\label{general_cores}
Let $\mathbb{Q}$ denote the rational numbers with the standard order. $\mathbb{Q}$ is an $\mathcal{I}$-core. On the other hand, $0:\mathbb{Q}\to\{0\}$ is a comparative retraction.

For another example, consider the space $X=\{(a,b):a\in\{0,1\},b\in\mathbb{N}\}$ with order given by: $(a_1,b_1)\leq (a_2,b_2)$ iff $(a_1,b_1)=(a_2,b_2)$ or $b_1<b_2$. $X$ is an $\mathcal{I}$-core. But consider the retractions $f_1,f_2$ given by $f_1(a,b)=(0,a+b)$ ($+$ indicates the standard addition of natural numbers), $f_2(a,b)=(0,1)$ for all $(a,b)\in X$. We have: $\id_{X}\leq f_1$, so $X$ is not a $\mathcal{C}$-core. Moreover, $f_1 \geq f_2$, so $X$ is contractible.
\end{example}

Since in this article we are mostly concerned with fp-spaces and locally finite spaces, the following proposition guarantees that in our setting $\mathcal{I}$-cores and $\mathcal{C}$-cores are the same. The proposition was proved in \cite{stong} for finite spaces.

\begin{proposition}\label{finite-chain_cores}
Let $X$ be a finite-chains $\mathcal{I}$-core (with the distinguished point $p$). If $f\sim\operatorname{id}_X$ for some $f\in C(X,X)$, then $f=\operatorname{id}_X$. Therefore, a finite-chains $\mathcal{I}$-core is a $\mathcal{C}$-core.
\end{proposition}
\begin{proof}
Let $f:X\to X$ (or $f:(X,p)\to (X,p)$) be a continuous map. Suppose $f\leq \id_{X}$. We will show $f=\id_X$.

Since $X$ satisfies DCC and Lemma \ref{max_is_acc} holds, we may use Noetherian induction. Take an $y\in X$ and assume that $f(x)=x$ for all $x<y$. By continuity of $f$, $x=f(x)\leq f(y)$ for $x<y$. But $f\leq\id_X$, so $f(y)\leq y$. If $f(y)<y$, then $y$ must be a down-beat point over $d_y=f(y)$. Otherwise, an $x<y$ would exist with $x\not\leq f(y)$ and $f$ would not be order-preserving. But $X$ is an $\mathcal{I}$-core, so it contains no beat points. Therefore, $f(y)=y$. By induction, $f=\id_X$. 

If $f\geq \id_X$, then $f=\id_X$ by the same reasoning as above. Therefore, $X$ is an $\mathcal{U}\cup\mathcal{D}$-core and thus a $\mathcal{C}$-core.
\end{proof}

In \cite{stong} Stong shows that for any finite space $X$ it is possible to reduce it to its $\mathcal{I}$-core $X^C$ by removing, one by one, beat points. Since the used retractions are comparative and Lemma \ref{finite_homotopies_lemma} holds, the resulting core is a strong deformation retract of the given finite space (and thus the space and its core are homotopy equivalent). Moreover, because of Proposition \ref{finite-chain_cores}, any function $f\in C(X^C,X^C)$ comparable to $\operatorname{id}_{X^C}$ is equal to $\operatorname{id}_{X^C}$. But $X^C$ is finite, so $C(X^C,X^C)$ is Alexandroff and that means $\operatorname{id}_{X^C}$ is an isolated point in $C(X^C,X^C)$. In particular, any $f\in C(X^C,X^C)$ homotopic to $\operatorname{id}_{X^C}$ is equal to $\operatorname{id}_{X^C}$. Now it easily follows (like in the proof of Corollary \ref{homotop_homeo}) that any two finite spaces are homotopy equivalent if and only if their $\mathcal{I}$-cores are homeomorphic. (Using the same arguments one can show that finite spaces with distinguished points $(X,p),(Y,q)$ are homotopy equivalent if and only if their cores $(X^C,p), (Y^C,q)$ are homeomorphic.)

From \cite{farley} we cite (in a version modified to our needs) the Li-Milner theorem, which allows us to generalise the above results to chain-complete posets without infinite antichains.
\begin{theorem}[cf. Theorem 6.11 of \cite{farley}]\label{li-milner}
For every chain-complete partially ordered set $X$ with no infinite antichains exists a finite $\mathcal{C}$-core $X^C\subseteq X$ that is a strong deformation retract of $X$. (In fact, $X$ is $\mathcal{C}$-dismantlable to $X^C$, in the sense of Definition \ref{def-dismantling}, in finitely many steps.)
\end{theorem}
Therefore, any two chain-complete ordered sets without infinite antichains are homotopy equivalent iff their cores (which exist and are finite) are homeomorphic.

\begin{question}
Is there a version of the above theorem for posets with distinguished points?
\end{question}

Now we will follow Stong's path in the setting of fp-spaces and bp-spaces. We shall describe a process which we call \cC-dismantling. It is a special case of the following, more general definition of dismantling a poset.

\begin{definition}[cf. Exercise 24 in Chapter 4 of \cite{schroder}]\label{def-dismantling}
Let $\mathcal{\gamma}$ be an ordinal and $X$ be an Alexandroff space (with the distinguished point $p$). Let $\{r_\alpha:X_\alpha\to X_{\alpha+1}\}_{\alpha<\gamma}$ be a family of retractions from a class $\mathcal{R}$ (keeping the distinguished point fixed) such that $X_0=X$, $X_{\alpha+1}=r_{\alpha}(X_{\alpha})$ for all $\alpha<\gamma$ and $X_{\alpha}=\bigcap_{\beta<\alpha}X_{\beta}$ for limit ordinals $\alpha<\gamma$. By transfinite induction we define a family of retractions $\{R_\alpha:X\to X_\alpha\}_{\alpha\leq\gamma}$:
\begin{enumerate}
\item $R_0=\id_X$,
\item $R_{\alpha+1}=r_{\alpha}\circ R_{\alpha}$,
\item for a limit ordinal $\alpha$ and an $x\in X$, if there exists $\beta_0<\alpha$ such that $R_{\beta}(x)=R_{\beta_{0}}(x)$ for all $\beta_{0}\leq\beta<\alpha$, then $R_{\alpha}(x)=R_{\beta_{0}}(x)$.
\end{enumerate}
The induction ends when $R_\gamma$ is defined or when $R_\alpha$ cannot be fully defined for some limit ordinal $\alpha$. In the first case we say the family $\{r_\alpha\}_{\alpha<\gamma}$ is \textit{infinitely composable} and $X$ is \textit{$\mathcal{R}$-dismantlable to $X_\gamma$ (in $\gamma$ steps)}. In the second case we say the family $\{r_\alpha\}_{\alpha<\gamma}$ is \textit{not infinitely composable} and if the induction stopped at $\gamma$ and $X_\gamma=\emptyset$, we say $X$ is \textit{$\mathcal{R}$-dismantlable to the empty set (in $\gamma$ steps)}.
\end{definition}

\begin{remark}\label{limits_remark}
By the standard construction of categorical limits, it is clear that the set $X_{\gamma}=\bigcap_{\alpha<\gamma}X_{\alpha}$ together with the family of inclusions $\{J_{\beta}:X_{\gamma}\hookrightarrow X_\alpha\}_{\alpha<\gamma}$ is the limit of the inverse system $\{i_{\alpha,\beta}:X_\beta\hookrightarrow X_\alpha\}_{\alpha<\beta<\gamma}$ of inclusions.

Suppose the sequence $\{r_\alpha\}_{\alpha<\gamma}$ is infinitely composable. We show $X_{\gamma}$ together with the family $\{R_{\alpha,\gamma}=R_{\gamma}|_{X_{\alpha}}:X_\alpha\to X_\gamma\}_{\alpha<\gamma}$ of retractions is also the colimit of the direct system $\{r_{\alpha,\beta}:X_\alpha\to X_\beta\}_{\alpha<\beta<\gamma}$ of retractions induced in the obvious way by the sequence $\{r_\alpha\}_{\alpha<\gamma}$.

It is trivial that $R_{\beta,\gamma} \circ r_{\alpha,\beta}=R_{\alpha,\gamma}$ for all $\alpha<\beta<\gamma$. If $Y, \{S_\alpha:X_\alpha\to Y\}_{\alpha<\gamma}$ are such that $S_{\beta} \circ r_{\alpha,\beta} = S_{\alpha}$ for all $\alpha<\beta<\gamma$, then the map $F:X_{\gamma} \to Y$ given by $F=S_{0} \circ J_{0}$ has the property that $F \circ R_{\alpha, \gamma} = S_{\alpha}$ for all $\alpha<\gamma$. If $F':X_{\gamma} \to Y$ is another map with this property, then $S_{0} = F' \circ R_{0, \gamma}$ and thus $S_{0} \circ J_{0} = F' \circ R_{0, \gamma} \circ J_{0} = F' \circ \id_{X_{\gamma}} = F'$. 

If the sequence $\{r_\alpha\}_{\alpha<\gamma}$ is not infinitely composable, but it $\mathcal{C}$-dismantles $X$ to the empty set, then the colimit of $\{r_{\alpha,\beta}:X_\alpha\to X_\beta\}_{\alpha<\beta<\gamma}$ is (by the usual construction of colimits) a one point space.
\end{remark}

\begin{definition}
An Alexandroff space $X$ (with distinguished point $p$) is \textit{countably $\mathcal{C}$-dismantlable} (or, for short, \textit{$c\mathcal{C}$-dismantlable}) to $X'\subseteq X$ ($(X',p)\subseteq (X,p)$) if it is $\mathcal{C}$-dismantlable to $X'$ ($(X',p)$) in a countable number of steps.
\end{definition}

By Theorem \ref{countable_homotopies} and Remark \ref{retractions_remark}, if $X$ (or $(X,p)$) is \cC-dismantlable to $\emptyset\not=X'\subseteq X$ (or $(X',p)\subseteq (X,p)$), then $X'$ (or $(X',p)$) is a strong deformation retract of $X$. (From the statement of the theorem it only follows that it is a deformation retract. However, if one examines the construction of the homotopy in the proof of the theorem, it becomes clear the word "strong" may be added.)

There are different ways one can try to \cC-dismantle a poset. The one we shall now present follows the idea of removing beat points. However, contrary to the finite case, we remove many of these at once. If the setting is right, this allows the process to finish after a relatively small number of steps. (One another idea of how to dismantle a poset is that of perfect sequences and retractions onto ``good'' subsets, cf. the proof of the Li-Milner theorem in \cite{farley} or \cite{schroder}.) 

\begin{definition}
Let $X$ be a finite-chains space (with distinguished point $p$). Let $u_X:X\to X$ be given by: 
\[u_X(x)=\begin{cases}u_x & \text{for } x \text{ an up beat point under } u_x \\ x & \text{otherwise} \end{cases}.\]
Since $u_X(x)\geq x$ for every $x\in X$ and $X$ is a finite-chains space, for every $x\in X$ exists an $N_x\in\mathbb{N}$ such that $(u_X)^n(x)=(u_X)^{N_x}(x)$ for every $n\geq N_x$. Let $U_X:X\to U_X(X)$ be given by $U_X(x)=(u_X)^{N_x}(x)$. It is straightforward to check $U_X$ is an up-retraction.

Dually we define the down retraction $D_X:X\to D_X(X)$.
\end{definition}

\begin{definition}
For an ordinal $\gamma$ and a finite-chains space $X$ we define by transfinite induction a sequence of retractions $\{r_{\alpha}:X_{\alpha}\to X_{\alpha+1}\}_{\alpha<\gamma}$, which we will call \textit{the standard sequence of $X$ (of length $\gamma$)}.

Let $X_0=X$,  $X_{\alpha+1}=r_{\alpha}(X)$ and $X_{\alpha}=\bigcap_{\alpha<\beta}X_{\beta}$ if $\alpha$ is a limit ordinal. For $\alpha=0$ or $\alpha$ a limit ordinal and $n$ a finite ordinal let 
\[r_{\alpha+n}=\begin{cases}D_{X_{\alpha+n}} & \text{if } n \text{ is even}\\ U_{X_{\alpha+n}} & \text{if } n \text{ is odd}\\ \end{cases}.\]
\end{definition}

The following theorem is the first of the two main results of this section.
\begin{theorem}
If $X$ is an fp-space (with the distinguished point $p$), then $X$ may be $\mathcal{C}$-dismantled to a $\mathcal{C}$-core $X^C\subseteq X$. If $X$ is in addition countable, then the dismantling may be chosen so that it takes countably many steps.

If $X$ is a bp-space (with the distinguished point $p$), then $X$ may be $\mathcal{C}$-dismantled to a $\mathcal{C}$-core $X^C\subseteq X$ in finitely many steps.
\end{theorem}
\begin{proof}
We may assume $X$ is infinite. Let $\Omega$ be the first ordinal of cardinality greater than $X$. Let $\{r_{\alpha}:X_{\alpha}\to X_{\alpha+1}\}_{\alpha<\Omega}$ be the standard sequence of $X$ of length $\Omega$.

First we show that $\{r_\alpha\}$ is infinitely composable. Indeed, it is easy to see that if for some limit ordinal $\gamma$ and $x\in X$ no ordinal $\beta_0$ existed such that $R_\beta(x)=R_{\beta_0}(x)$ for all $\beta_0<\beta<\gamma$, then for every $\alpha_0<\gamma$ the sequence $\{R_\alpha(x)\}_{\alpha_0<\alpha<\gamma}$ would contain an infinite s-path.

Now we show that for $X$ an fp-space the standard sequence is constant beginning with some $\alpha_0<\Omega$. That means $X_{\alpha_0}$ is a core and thus $X$ is $\mathcal{C}$-dismantlable to $X_{\alpha_0}$ in $\alpha_0$ steps. If $X$ is countable, then $\Omega=\omega_1$, the first uncountable cardinal, and thus $\alpha_0$ is countable.

It is easy to see that if for some ordinal $\alpha<\Omega$ both maps $r_{\alpha}, r_{\alpha+1}$ are identities, then the sequence is constant begining with $\alpha$. Therefore, if for no $\alpha_0<\Omega$ the sequence was constant beginning with $\alpha_0$, then for every ordinal $\alpha<\Omega$ we could choose an $x_{\alpha}\in X$ such that $x_{\alpha}$ was removed by one of the maps $r_{\alpha+i}, i=0,1$. It is now easy to see that the set $\{x_\alpha:\alpha<\Omega\}\subseteq X$ would be of the same cardinality as $\Omega$ (which is greater than the cardinality of $X$) -- a contradiction.

If $X$ is a bp-space with length of s-paths bounded by some $n\in\mathbb{N}$, then we may show the standard sequence is constant beginning with $2n+2$. Suppose it is not so. Thus there exists an $m\geq 2n+2$ with $r_m\not=\id_{X_m}$. This means $r_m(x_0)\not=x_0$ for some $x_0\in X$. By the definition of the standard sequence, $x_0$ is a beat point. Since $x_0$ was not removed by any of the previous retractions, there exists an $x_1\in X_{m-1}\cup X_{m-2}$ with $x_1\sim x_0$. The same reasoning lets us define inductively $x_0\sim x_1\sim x_2\sim \ldots \sim x_{n}$. But that is an s-path of length greater than $n$.
\end{proof}

One can easily find examples of uncountable fp-spaces that are not bp-spaces and nonetheless are \cC-dismantlable to a core. As the following Example \ref{monstruous_spider} shows, fp-spaces exist that are not \cC-dismantlable to a core and nonetheless have a $\mathcal{C}$-core as their strong deformation retract. We may ask whether this holds in the general.

\begin{question}
Is it true that for every fp-space $X$ exists a $\mathcal{C}$-core that is a strong deformation retract of $X$? More generally, if $X$ is a space $\mathcal{C}$-dismantlable to a $\mathcal{C}$-core, is the core necesarrily a strong deformation retract of $X$?
\end{question}

\begin{example}[cf. \cite{rutkowski}]\label{monstruous_spider}
We will define by transfinite induction an fp-space $S$ which cannot be \cC-dismantled to a core and a homotopy $H:S\times I\to S$ from $\id_S$ to a constant map.

Let $S_0=\{s_0\}$ and $H_0(s_0,t)=s_0$ for all $t\in I$.

If $\alpha=0$ or $\alpha$ is a limit ordinal and $n\geq 0$ is a finite ordinal, then $S_{\alpha+n+1}=S_{\alpha+n}\cup\{s_{\alpha+n+1}\}$. If $n$ is even, we add the comparability $s_{\alpha+n+1}>s_{\alpha+n}$ and put \[H_{\alpha+n+1}(x,t)=\begin{cases}H_{\alpha+n}(x,2t)&\text{for } x\not=s_{\alpha+n+1}\text{ and } t\leq\frac{1}{2}\\ s_{\alpha+n} & \text{for } x\not=s_{\alpha+n+1}\text{ and } t\in(\frac{1}{2},1)\\s_{\alpha+n+1} & \text{otherwise} \end{cases}.\] If $n$ is odd, we add the comparability $s_{\alpha+n+1}<s_{\alpha+n}$ and put \[H_{\alpha+n+1}(x,t)=\begin{cases}H_{\alpha+n}(x,2t)&\text{for } x\not=s_{\alpha+n+1}\text{ and } t\leq\frac{1}{2}\\ s_{\alpha+n+1} & \text{otherwise}\end{cases}.\]

If $\alpha$ is a limit ordinal, then let $S_{\alpha}=\{s_{\alpha}\}\cup\bigcup_{\beta<\alpha}S_{\beta}\times\{\beta\}$ with the comparabilities $s_{\alpha}<(s_{\beta},\beta)$ for all $\beta<\alpha$ added. We define
\[H_{\alpha}(x,t)=\begin{cases}(H_{\beta}(y,2t),\beta) & \text{for } x=(y,\beta), \beta<\alpha \text{ and } t\leq\frac{1}{2}\\ s_{\alpha} & \text{otherwise}\end{cases}.\]

By transfinite induction we prove $S_{\alpha}$ is an fp-space for every ordinal $\alpha$. For $\alpha=0$ this is trivial. If $\alpha=\beta+1$ and the statement holds for $\beta$, then it also holds for $\alpha$, since any s-path in $S_{\alpha}$ is either an s-path in $S_{\beta}$ or can be obtained from an s-path in $S_{\beta}$ by adding $s_{\alpha}$ at the beginning or at the end of the s-path. If $\alpha$ is a limit ordinal and the statement is true for all $\beta<\alpha$, then it also holds for $\alpha$, since an s-path in $S_{\alpha}$ intersects at most two of the sets $S_{\beta}\times\{\beta\}, \beta<\alpha$ (and intersection with each of these sets is also an s-path), so its length is at most the sum of the lengths of intersections enlarged by 1. 

Let $\alpha$ be an uncountable ordinal. We put $S=S_\alpha$, $H=H_\alpha$. $S$ cannot be \cC-dismantled to a core. This can be checked by examining functions comparable to the identity on $S$, then functions comparable to these functions, and so on. The shortest possible dismantling has uncountably many steps. Nevertheless, $H$ is a homotopy from $\id_{S}$ to the constant map $s_{\alpha}$.

Another idea of how to obtain such a homotopy is to apply Theorem \ref{one_path_contractible}. 
\end{example}

We now proceed to the second main theorem.

\begin{theorem}\label{homotopic_are_id}
If $X$ is an fp-space (with the distinguished point $p$) that is a $\mathcal{C}$-core, then the connected component of $\id_X$ in $C(X,X)$ is a singleton.
\end{theorem}
\begin{proof}
We will first show that for every $x\in X$ exists a finite set $\mathbf{A}_x\subseteq X$ with $x\in \mathbf{A}_x$ that has the following property: if $y\in \mathbf{A}_x$, then either $y=p$, or $|\mathbf{A}_x\cap\max(y\!\downarrow\smallsetminus\{y\})|\geq 2$ if $y$ is not minimal in $X$ and $|\mathbf{A}_x\cap\min(y\!\uparrow\smallsetminus\{y\})|\geq 2$ if $y$ is not maximal in $X$.

If $x=p$, then put $\mathbf{A}_x=\{p\}$. Otherwise, label $x=()$ (where $()$ is the empty sequence) and let $A_0=\{()\}$. Suppose we have defined finite sets $A_m\subseteq X$ for all $0\leq m<n$ and the elements of $A_m$ are labeled by sequences of length $m$ with terms from the set $\{1,2,3,4\}$. In particular, $A_{n-1}=\{x_1,x_2,\ldots,x_k\}$ is defined and every element $x_i\in A_{n-1}$ has the label $(x_i^1,x_i^2,\ldots,x_i^{n-1})$. For $i=1,\ldots,k$ we will define sets $A_{n}^i$. If $x_i\not=p$ and $x_i$ is not maximal, then there exist two distinct elements in $\min(x\!\uparrow\smallsetminus\{x_i\})$. Label these elements $(x_i^1,x_i^2,\ldots,x_i^{n-1},1)$ and $(x_i^1,x_i^2,\ldots,x_i^{n-1},2)$. If $x_i\not=p$ and $x_i$ is not minimal, then there exist two distinct elements in $\max(x_i\!\downarrow\smallsetminus\{x_i\})$. Label the elements $(x_i^1,x_i^2,\ldots,x_i^{n-1},3)$ and $(x_i^1,x_i^2,\ldots,x_i^{n-1},4)$. Let the set $A_{n}^i$ consist of those of the four newly labelled elements that do not belong to $(\bigcup_{r\leq n-1}A_{r}\cup\bigcup_{j<i}A_{n}^j)$. Now let $A_n=\bigcup_{i=1}^{k}A_n^i$ and $\mathbf{A}_x=\bigcup_{n=0}^{\infty}A_n$.

We must only show that $\mathbf{A}_x$ is finite, other properties that we require from $\mathbf{A}_x$ are clearly satisfied. Suppose it is infinite. Then there exists an $n_1\in\{1,2,3,4\}$ such that the set of elements of $\mathbf{A}_x$ with labels beginning with $n_1$ is infinite. If there exist $n_1,n_2,\ldots,n_k\in\{1,2,3,4\}$ such that the set of elements of $\mathbf{A}_x$ with labels beginning with the sequence $(n_1,n_2,\ldots,n_k)$ is infinite, then there exists an $n_{k+1}\in\{1,2,3,4\}$ such that the set of elements of $\mathbf{A}_x$ with labels beginning with the sequence $(n_1,n_2,\ldots,n_k,n_{k+1})$ is infinite. By induction we can define an infinite s-path $( (), (n_1), (n_1,n_2), (n_1,n_2,n_3), \ldots)$ in $X$ -- a contradiction. Thus, $\mathbf{A}_x$ must be finite.

For every $x\in X$ the set $\bigcap_{y\in \mathbf{A}_x}[y,y\!\downarrow]$ is an open neighbourhood of $\id_X$. We will show that the neighbourhood is also closed. Moreover, if $f\in \bigcap_{y\in \mathbf{A}_x}[y,y\!\downarrow]$, then $f|_{\mathbf{A}_x}=\id_{\mathbf{A}_x}$. Since the connected component of $\id_X$ is a subset of the quasi-component of $\id_X$ (= the intersection of all clopen sets containing $\id_X$) and from the above it follows that the quasi-component is contained in $\bigcap_{x\in X}\bigcap_{y\in \mathbf{A}_x}[y,y\!\downarrow]=\{\id_X\}$, this will finish the proof.

Let $x\in X$ and let $f\in \bigcap_{y\in \mathbf{A}_x}[y,y\!\downarrow]$. Suppose there is some $a_0\in \mathbf{A}_x$ with $f(a_0)<a_0$. If we have found $(a_n)_{n=1}^{m}\subseteq\mathbf{A}_x$ with the property that $f(a_n)<a_n$ for $0\leq n\leq m$ and $a_{k+1}<a_{k}$ for $0\leq k < m$, then, since $f(a_m)<a_m$, $a_m$ is not minimal. Therefore, there exists an $a_{m+1}\in \mathbf{A}_x$ with $a_{m+1}<a_{m}$ and $a_{m+1}\not\leq f(a_m)$. $f$ is order-preserving, so $f(a_{m+1})\leq f(a_m)$. Because $f\in [a_{m+1},a_{m+1}\!\downarrow]$, $f(a_{m+1})<a_{m+1}$. By induction we have defined a sequence $(a_n)_{n\in\mathbb{N}}\subseteq \mathbf{A}_x$ with $a_{n+1}<a_{n}$ -- a contradiction, since $\mathbf{A}_x$ satifies the DCC. So $f(a)=a$ for all $a\in \mathbf{A}_x$.

For every $x\in X$ and $g\in C(X,X)$ such that $g\not\in \bigcap_{y\in \mathbf{A}_x}[y,y\!\downarrow]$ we will now find an open neighbourhood $U$ of $g$ with $U\cap \bigcap_{y\in \mathbf{A}_x}[y,y\!\downarrow]=\emptyset$. (This will mean that $X\smallsetminus \bigcap_{y\in \mathbf{A}_x}[y,y\!\downarrow]$ is open and thus $\bigcap_{y\in \mathbf{A}_x}[y,y\!\downarrow]$ is closed.) Fix some $x$ and $g$ as above. By what we have shown above, $g\not\in \bigcap_{y\in \mathbf{A}_x}[y,y\!\downarrow]$ iff there exists an $a\in \mathbf{A}_x$ with $g(a)\not= a$. If $g(a)<a$ or $g(a)\not\sim a$, then we may take $U=[a,g(a)\!\downarrow]$. If $g(a)>a$, then we show there exists an $\overline{a}\in \mathbf{A}_x$ with $g(\overline{a})\not\sim\overline{a}$ and take $U=[\overline{a},g(\overline{a})\!\downarrow]$. Indeed, suppose such an $\overline{a}$ does not exist. Let $a_0=a$. If $(a_n)_{n=1}^{m}\subseteq\mathbf{A}_x$ are defined, $g(a_n)>a_n$ for $0\leq n\leq m$ and $a_{k+1}> a_k$ for $0\leq k < m$, then $a_m\in \mathbf{A}_x$ is not maximal. Therefore, there exists an $a_{m+1}\in \mathbf{A}_x$ with $a_{m+1}\not\geq g(a_m)$. Since $g(a_{m+1})\sim a_{m+1}$, we must have $g(a_{m+1})>a_{m+1}$. By induction we have defined an infinite, strictly increasing sequence -- a contradiction.
\end{proof}

The following corollary is now proved as in the finite case.
\begin{corollary}\label{homotop_homeo}
If $X,Y$ are fp-spaces (with distinguished points $p,q$) such that the $\mathcal{C}$-cores $X^C\subseteq X,Y^C\subseteq Y$ are their strong deformation retracts, then $X$ (or $(X,p)$) is homotopy equivalent to $Y$ (or $(Y,q)$) if and only if $X^C$ is homeomorphic to $Y^C$ ($(X^C,p)$ is homeomorphic to $(Y^C,q)$). 
\end{corollary}
\begin{proof}
We know $X, Y$ are homotopy equivalent to their cores. Therefore, if the cores are homeomorphic, they are homotopy equivalent, and in turn the spaces $X, Y$ are homotopy equivalent. For the reverse implication, suppose $X$ and $Y$ are homotopy equivalent. Then their cores are homotopy equivalent, which means that there exist maps $f:X^C\to Y^C$ and $g:Y^C\to X^C$ with $f\circ g\simeq \id_{Y^C}$, $g\circ f\simeq\id_{X^C}$ (relative to base points). But from Theorem \ref{homotopic_are_id} it follows that $f\circ g=\id_{Y^C}$ and $g\circ f=\id_{X^C}$, so we have the desired homeomorphism.
\end{proof}

\begin{remark}\label{fpp_hom_invariant}
From Theorem \ref{li-milner} and Exercise 24 a) in Chapter 4 of \cite{schroder} it follows that the fixed point property is a homotopy invariant in the category of chain-complete posets without infinite antichains. The same Exercise 24 and Corollary \ref{homotop_homeo} show the same is true for $T_0$ bp-spaces and countable fp-spaces. From Example \ref{khalimsky_contractible} and Corollary \ref{loc_fin_non_fpp} we conclude, though, that the fixed point property is not a homotopy invariant in the category of locally finite $T_0$ spaces.
\end{remark}

What about cores in arbitrary posets? As Example \ref{general_cores} shows, $\mathcal{I}$-cores do not seem to be of much use in investigating homotopy types of Alexandroff spaces in general. $\mathcal{C}$-cores may be a more useful tool, although even quite simple $\mathcal{C}$-cores exist that do not have the property that their only self-map homotopic to the identity is the identity itself. As we have seen in Example \ref{khalimsky_contractible}, the Khalimsky line is one such space. This shows that Theorems 3.5 and 3.6 of \cite{arenas} are not true, and in turn the proof of Corollary 3.7 of \cite{arenas} is not valid. But even in the setting of all Alexandroff spaces it seems quite possible for the conclusion of Corollary 3.7 to be true. So we pose the following question.
\begin{question}
Is it true that an Alexandroff space $X$ is contractible if and only if some point of $X$ is a strong deformation retract of $X$? In particular, is it so for locally finite $X$? (By Corollary \ref{when_height_one_is_contractible}, the answer is yes for spaces of height 1.)
\end{question}

Even locally finite spaces exist that cannot be $\mathcal{C}$-dismantled to a core. This is the case for the Khalimsky half-line. The half-line is, however, \cC-dismantlable to the empty set. More generally, for connected, locally finite spaces one has the following proposition.
\begin{proposition}\label{dism_core_or_empty}
If $X$ is a connected, locally finite Alexandroff space, then $X$ is $\mathcal{C}$-dismantlable in $\omega$ steps either to a $\mathcal{C}$-core or to the empty set. If $X$ is a connected, locally finite Alexandroff space with the distinguished point $p$, then $X$ is $\mathcal{C}$-dismantlable to a $\mathcal{C}$-core in $\omega$ steps.
\end{proposition}
\begin{proof}
Let $\{r_{\alpha}:X_{\alpha}\to X_{\alpha+1}\}_{\alpha<\omega}$ be the standard sequence of $X$ of length $\omega$. 

First suppose $\{r_\alpha\}$ is infinitely composable. If $X_{\omega}$ is not a $\mathcal{C}$-core, then there exists a beat point $x\in X_{\omega}$. Without lack of generality, assume $x$ is an up-beat point. Since $x$ has not been removed by any of the retractions $r_n, n<\omega$, then for every $n<\omega$ the set $\min((x\!\uparrow\smallsetminus\{x\})\cap X_n)$ has more than one element. On the other hand, $x\!\uparrow$ is finite. Thus, if for $y\in (x\!\uparrow)$ we define $n_y$ to be the ordinal such that $y$ was removed by $r_{n_y}$, or $n_y=0$ if $y\in X_{\omega}$, then $N=\max\{n_y:y\in x\!\uparrow\}$ exists and is finite. Therefore, $\min((x\!\uparrow\smallsetminus\{x\})\cap X_{\omega})=\min((x\!\uparrow\smallsetminus\{x\})\cap X_{m})$ is a singleton for all $m>N$ -- a contradiction.

Now suppose $\{r_\alpha\}$ is not infinitely composable. We show that $X_{\omega}=\emptyset$. This is the case if and only if for every $y\in X$ and every $n<\omega$ exists $n<m<\omega$ such that $R_m(x)\not=R_n(x)$. Since $\{r_\alpha\}$ is not infinitely composable, this happens for at least one $x_0\in X$. Let $x_n=R_n(x_0)$ for all $n<\omega$. Let $y\in B(x_0,1)$; without lack of generality $y>x_0$. We will define by induction a sequence $S(y)=(y^n)_{n\in\mathbb{N}}$ with $y^n\not=y^{n+1}$, $y^n=R_{k^n}(y^0)$ for some $k^n<\omega$ and $k^{n+1}>k^{n}$ for all $n\in\mathbb{N}$. Let $y^0=y$, $k^0=0$. Suppose for all $n\leq m$ elements $y^{n}=R_{k^n}(y^0)$ have been defined with the property that $y^n\geq x_{k^n}$. Since $y^{m}\!\downarrow$ is finite, there is some $k^{m+1}<\omega$ with $x_{k^{m+1}}\not\in y^{m}\!\downarrow$. Therefore, for $R_{k^{m+1}}$ to be continuous, \[y^{m+1}=R_{k^{m+1}}(y^0)\not=y^{m}=R_{k^{m}}(y^0).\] The sequence $S(y)$ is defined. By the fact that $X=\bigcup_{n\in\mathbb{N}}B(x_0,n)$ and simple induction ,we conclude that for every $y\in X$ exists an infinite sequence $S(y)$, which ends the proof for spaces without distinguished point. 

If $p$ is a distinguished point in $X$, then for every $n$ we have $R_n(p)=p$. Therefore, by the above reasoning, $\{r_\alpha\}$ must be infinitely composable.
\end{proof}

The following proposition shows that dismantlability to the empty set may be an useful notion.
\begin{proposition}\label{height_one_empty_set}
If $X$ is a connected Alexandroff space of height 1 that is $\mathcal{C}$-dismantlable to the empty set, then $X$ is contractible.
\end{proposition}
\begin{proof}
We will show that $X$ contains no crowns. Then the result follows from Corollary \ref{when_height_one_is_contractible}. Let $\{r_\alpha\}_{\alpha<\gamma}$ be a fixed dismantling of $X$ to the empty set. Suppose $\{x_0,x_1,\ldots,x_{m-1}\}\subseteq X$ is a crown. In the following, by $x_{i+k}$ we mean $x_{i+k\operatorname{mod} m}$. Define $n_i=\min\{\alpha<\gamma:r_{\alpha}(x_i)\not=x_i\}$ for $0\leq i< m$ and let $i_0$ be such that $n_{i_0}=\min(n_0,\ldots,n_m)$. We may assume $x_{i_0}\in\min(X)$ (for $x_{i_0}\in\max(X)$ the proof is analogous). Then $r_{n_{i_0}}(x_{i_0})=y$ for some $y> x_{i_0}$. Since $x_{i_0+1}, x_{i_0-1}\geq x_{i_0}$, $r_{n_{i_0}}(x_{i_0+1})\geq y \leq r_{n_{i_0}}(x_{i_0-1})$. But $y\in\max(X)$, so $r_{n_{i_0}}(x_{i_0+1})= y = r_{n_{i_0}}(x_{i_0-1})$. Thus, $y\sim x_{i_0+1}, x_{i_0-1}$, but $y\not\geq x_{i_0+1}, x_{i_0-1}$, since $x_{i_0+1}, x_{i_0-1}$ are two distinct, maximal elements. Therefore, $y<x_{i_0+1}, x_{i_0-1}$, which contradicts $y\in\max(X)$. 
\end{proof}

Proposition \ref{height_one_empty_set}  and Remark \ref{limits_remark} motivate the following question.
\begin{question}
If a connected Alexandroff space $X$ is $\mathcal{C}$-dismantlable to the empty set, is $X$ contractible? What if $X$ is locally finite or the dismantling process takes countably many (or $\omega$) steps?
\end{question}

We close this section with a simple application of a fixed-point-theoretic result about cores and dismantlings.
\begin{proposition}
If $X$ is an infinite, locally finite Alexandroff space, then $X$ is not $\mathcal{C}$-dismantlable to a finite $\mathcal{C}$-core.
\end{proposition}
\begin{proof}
We may assume $X$ is connected. Suppose $X$ is $\mathcal{C}$-dismantlable to a finite $\mathcal{C}$-core $X^C$ and the dismantling is given by the sequence of retractions $\{r_\alpha:X_\alpha\to X_{\alpha+1}\}_{\alpha<\gamma}$. Let $X'$ be the set $X\cup \{*\}$ with order induced from $X$ and comparabilities $x<*$ added for all $x\in\bigcup_{y\in X^C}y\!\downarrow$. It's easy to see $X'$ is locally finite, and thus, by Corollary \ref{loc_fin_non_fpp}, it does not have the fixed point property. The sequence $\{r_\alpha ':X_\alpha\cup\{*\}\to X_{\alpha+1}\cup\{*\}\}_{\alpha<\gamma}$ given by $r_\alpha ' (x)=\begin{cases}r(x) & \text{for } x\not=*\\ * & \text{for } x=*\end{cases}$ gives a $\mathcal{C}$-dismantling of $X'$ to $X^C\cup \{*\}$. But $X^C\cup \{*\}$ has a greatest element and thus has the fixed point property. It follows from Exercise 24 a) in Chapter 4 of \cite{schroder} that $X'$ also has the fixed point property. A contradiction.
\end{proof}
\end{section}

%
%

\begin{section}{Other results}

\begin{subsubsection}*{H-spaces}
Results of §5 of \cite{stong} on finite H-spaces of type I (for terminology, see the cited work) hold for bp-spaces and countable fp-spaces with only a slight change in the proof of Proposition 13 (instead of using the argument that unions of the sets $D_r$, $D_r'$ are infinite, we use the fact that they contain infinite s-paths). The following theorem holds:

\begin{theorem}[cf. \cite{stong}]
Let $X$ be a bp-space or a countable fp-space. A necessary and sufficient condition that there exists an H-space of type I $(X,p,\nu)$ is that $p$ be a strong deformation retract of its component in $X$. A necessary and sufficient condition that there is a H-structure of type I on $X$ for some base point is that a component of $X$ is contractible.
\end{theorem}
\end{subsubsection}

\begin{subsubsection}*{Maps from polyhedrons to Alexandroff spaces}
All results of §6 of \cite{stong} hold (with proofs unaltered) when $F$ being a finite space is replaced with $F$ being an Alexandroff space (not necesarrily $T_0$). In particular, we have the following theorem.

\begin{theorem}[cf. \cite{stong}]
Let $K$ be a finite simplicial complex, $L$ a closed subcomplex, $X$ an Alexandroff space, $p\in X$ and $\mathfrak{A}$ the space of continuous functions $h:(|K|,|L|)\to (X,p)$ with the compact-open topology.
\begin{enumerate}
\item If $f\in\mathfrak{A}$, there exists a $g\in\mathfrak{A}$ such that $\{h\in\mathfrak{A}:h\leq g\}$ is a neighbourhood of $f$ in $\mathfrak{A}$.
\item If $f,f'$ are homotopic (relative to $L$), there exist elements $\phi_i\in\mathfrak{A}$, $0\leq i\leq s$ with $\phi_0=f$, $\phi_s=f'$ and $\phi_i\sim\phi_{i+1}$ for $0\leq i<s$.
\end{enumerate}
\end{theorem}
\end{subsubsection}

\begin{subsubsection}*{$\gamma$-points}
Many of the results concerning so-called $\gamma$-points which in the article \cite{barmak} by Barmak and Minian were stated for finite spaces also hold for Alexandroff spaces in general, with proofs unchanged. In particular, we have the following definition and theorems.

\begin{definition}[cf. \cite{barmak}]
A point $x$ of an Alexandroff $T_0$-space $X$ is a \textit{$\gamma$-point} if $\hat{C}_x=\{y\in X:y\sim X\}\smallsetminus\{x\}$ is homotopically trivial (i.e. if all its homotopy groups are trivial). 
\end{definition}

\begin{proposition}[cf. \cite{barmak}]
If $x\in X$ is a $\gamma$-point, the inclusion $i:X\smallsetminus \{x\} \hookrightarrow X$ is a weak homotopy equivalence.
\end{proposition}

\begin{proposition}[cf. \cite{barmak}]
Let $x$ be a point of an Alexandroff $T_0$-space $X$. The inclusion $i:X\smallsetminus \{x\} \hookrightarrow X$ induces isomorphisms in all homology groups if and only if the subspace $\hat{C}_x$ is acyclic.
\end{proposition}

\begin{theorem}[cf. \cite{barmak}]
Let $X$ be an Alexandroff $T_0$-space, and $x\in X$ a point which is neither maximal nor minimal and such that $X\smallsetminus \{x\} \hookrightarrow X$ is a weak homotopy equivalence. Then $x$ is a $\gamma$-point.
\end{theorem}

The same authors in \cite{barmak1} also studied simple homotopy types of finite spaces. It seems an interesting question how far their results may be applied in the infinite case.
\end{subsubsection}
\end{section}

%
%

\end{document}